\documentclass{article}
\usepackage{proof,logic,latexsym,amssymb,times,mathptm}

\title{Approximating Propositional Calculi by\\
Finite-valued Logics}
\date{}

\author{Matthias Baaz\thanks{Technische Universit\"at Wien, Institut f\"ur Algebra und Computermathematik E118.2, A-1040
Vienna, Austria. Email: baaz@logic.at} \and Richard
Zach\thanks{University of Calgary, Department of Philosophy Calgary,
Alberta T2N 1N4, Canada. Email: rzach@ucalgary.ca}}

\let\impl\supset
\def\val{{\rm val}}
\let\Box\square

\def\qedchar{\quad\vrule height7pt width5pt depth-1pt}
\def\I{{\bf I}}
\def\J{{\bf J}}

\def\cd#1{\lceil #1 \rceil}
\def\tbox{\widetilde{\vphantom{l}\Box}}
\def\tdiamond{\widetilde{\vphantom{l}\Diamond}}
\def\MC{{\it MC}}

\def\CL{{\bf CL}}
\def\LL{{\bf LL}}
\def\IPC{{\bf IPC}}
\def\A{{\bf A}}
\def\IPL{{\bf IPL}}
\def\Taut{{\rm Taut}}
\def\Thm{{\rm Thm}}
\def\Frm{{\rm Frm}}
\def\Var{\mathord{\rm Var}}
\def\M{{\bf M}}
\def\AL{{\bf L}}
\def\AK{{\bf C}}
\def\SF{{\bf S4}}
\def\G{{\bf G}}
\def\Gm{\imm{{\bf G}_m}}
\def\LA{\imm{{\cal L}}}
\def\K{\imm{{\cal K}}}
\def\MVL{{\rm MVL}}
\def\better{\mathrel{\lhd}}
\def\bettereq{\mathrel{\unlhd}}
\let\best\bigtriangleup
\let\worst\bigtriangledown
\def\cbi{\mathop{\tilde{\leftrightarrow}}}
\def\Nex{{\scriptstyle \bigcirc}}

\begin{document}
\bibliographystyle{abbrv}
\maketitle

\thispagestyle{empty}

\subsection*{\centering Abstract}
{\em Bernays introduced a method for proving underivability results in
propositional calculi~\AK{} by truth tables.  In general, this
motivates an investigations of how to find, given a propositional
logic, a finite-valued logic which has as few tautologies as possible,
but which has all the valid formulas of the given logic as
tautologies.  It is investigated how far this method can be carried
using (1)~one or (2)~an infinite sequence of finite-valued logics.  It
is shown that the best candidate matrices for (1) can be computed from
a calculus, and how sequences for~(2) can be found for certain classes
of logics (including, in particular, logics characterized by Kripke
semantics).}

\section{Introduction}

The question of what to do when face to face with a new logical
calculus is an age-old problem of mathematical logic.  One usually
has, at least at first, no semantics.  For example, intuitionistic
propositional logic was constructed by Heyting only as a calculus;
semantics for it were proposed much later.  Currently we face a
similar situation with Girard's linear logic.  The lack of semantical
methods makes it difficult to answer questions such as: Are statements
of a certain form (un)derivable? Are the axioms independent? Is the
calculus consistent?  For logics closed under substitution many-valued
methods have often proved valuable since they were first used for
proving underivabilities by Bernays \cite{Bernays:26} in 1926 (and
later by others, e.g., McKinsey and Wajsberg; see also
\cite[\S~25]{Rescher:69}).  For the above-mentioned underivability
question it is necessary to find many-valued matrices for which the
given calculus is sound.  If a formula is not a tautology under such a
matrix, it cannot be derivable in the calculus.  It is also necessary,
of course, that the matrix has as few tautologies as possible in order
to be useful.

Such ``optimal'' approximations of a given calculus may also have
applications in computer science. In the field of artificial
intelligence many new (propositional) logics have been introduced.
They are usually better suited to model the problems dealt with in AI
than traditional (classical, intuitionistic, or modal) logics, but
many have two significant drawbacks: First, they are either given
solely semantically or solely by a calculus.  For practical purposes,
a proof theory is necessary; otherwise computer representation of and
automated search for proofs/truths in these logics is not feasible.
Second, most of them are intractable, and hopelessly so, provided the
polynomial hierarchy does not collapse.  For instance, many
nonmonotonic formalisms have been shown to be hard for classes above
NP \cite{EiterGottlob:92}.  Although satisfiability in many-valued
propositional logics is (as in classical logic) NP-complete
\cite{Mundici:87}, this is still (probably) much better.

On the other hand, it is evident from the work of
Carnielli~\cite{Carnielli:87} and H\"ahnle~\cite{Hahnle:93} on
tableaux, and Rousseau, Takahashi, and Baaz et
al.~\cite{BaazFermZach:93} on sequents, that finite-valued logics are,
from the perspective of proof {\em and} model theory, very close to
classical logic.  Therefore, many-valued logic is a very suitable
candidate if one looks for approximations, in some sense, of given
complex logics.

What is needed are methods for obtaining finite-valued approximations
of the propositional logics at hand.  It turns out, however, that a
shift of emphasis is in order here.  While it is the {\em logic} we
are actually interested in, we always are given only a {\em
representation} of the logic.  Hence, we have to concentrate on
approximations of the representation, and not of the logic per se.

What is a representation of a logic?  The first type of representation
that comes to mind is a calculus.  Hilbert-type calculi are the
simplest conceptually and the oldest historically.  We will
investigate the relationship between such calculi on the one hand and
many-valued logics or recursive sequences of many-valued logics on the
other hand.  The latter notion has received considerable attention in
the literature in the form of the following two problems: Given a
calculus~\AK,
\begin{enumerate}
\item find a minimal (finite) {\em normal} matrix for~\AK{} (relevant for
non-derivability and independence proofs), and
\item find a sequence of finite-valued logics
whose intersection equals the theorems of~\AK, and its converse,
given a sequence of finite-valued logics, find a calculus
for its intersection (exemplified by Ja{\'s}kowski's
sequence for intuitionistic propositional calculus, and by
Dummett's extension axiomatizing the intersection of the sequence
of G\"odel logics, respectively).
\end{enumerate}
For (1), of course, the best case would be a finite-valued logic~\M{}
whose tautologies {\em coincide} with the theorems of~\AK.  \AK{} then
provides an axiomatization of~\M.  This of course is not always
possible, at least for {\em finite}-valued logics.  Lindenbaum
\cite[Satz~3]{LukasiewiczTarski:30} has shown that any logic (in our
sense, given by a set of rules and closed under substitution) can be
characterized by an {\em infinite}-valued logic.  For a discussion of
related questions see also Rescher~\cite[\S~24]{Rescher:69}.

In the following we will consider these questions in a general
setting.  Consider a propositional Hilbert-type calculus~\AK.  First
of all, an optimal (i.e., minimal under set inclusion of the
tautologies) $m$-valued logic for which~\AK{} satisfies reasonable
soundness properties can be computed.  We call such a logic {\em
normal} for~\AK.  The next question is, can we find an approximating
sequence of $m$-valued logics in the sense of~(2)?  It is shown that
this is impossible for undecidable calculi~\AK, and possible for all
decidable logics closed under substitution.  This leads us to the
investigation of the {\em many-valued closure} $\MC(\AK)$ of~\AK,
i.e., the set of formulas which are true in all covers of~\AK. In
other words, if some formula can be shown to be underivable in~\AK{}
by a Bernays-style many-valued argument, it is not in the many-valued
closure.  Using this concept we can classify calculi according to
their many-valued behavior, or according to how good they can be dealt
with by many-valued methods.  In the best case $\MC(\AK)$ equals the
theorems of~\AK{} (This can be the case only if \AK{} is decidable).
Otherwise $\MC(\AK)$ is a proper superset of the theorems of~\AK.

We show that $\MC(\AK)$ is decidable if \AK~is {\em analytic} (This
does not imply that \AK{} itself is decidable; e.g., cut-free
propositional linear logic is known to be undecidable).  Two
axiomatizations \AK~and~$\AK'$ of the same logic may have different
many-valued closures $\MC(\AK)$ and $\MC(\AK')$ while being
model-theoretically indistinguishable. Hence, the many-valued closure
can be used to distinguish between \AK~and~$\AK'$ with regard to their
proof-theoretic properties.

Finally, we investigate some of these questions for other
representations of logics, namely for decision procedures and finite
Kripke models.  In these cases approximating sequences of many-valued
logics whose intersection equals the given logics can always be given.

\section{Propositional Logics}

\begin{defn}
A {\em propositional language~$\cal L$} consists of the following:
\begin{enumerate}
\item propositional variables:  $X_0$, $X_1$, $X_2$, \dots, $X_j$,
\dots{} ($j \in \omega$)
\item propositional connectives of arity $n_j$:
$\Box_0^{n_0}$,~$\Box_1^{n_1}$, \dots,~$\Box_r^{n_r}$. If $n_j = 0$,
then $\Box_j$ is called a {\em propositional constant}.
\item Auxiliary symbols: $($, $)$, and $,$ (comma).
\end{enumerate}
\end{defn}

Formulas and subformulas are defined as usual.
We denote the set of formulas over a
language $\cal L$ by $\Frm({\cal L})$.
By $\Var(A)$ we mean the set of propositional
variables occurring in~$A$.

\begin{defn}\label{defn:proplogic}
A {\em propositional Hilbert-type calculus}~\AK{} in the language ${\cal L}$ is given by
\begin{enumerate}
\item A finite set $A(\AK) \subseteq \Frm(\LA)$ of axioms.
\item A finite set $R(\AK)$ of rules of the form
\[
\infer[r]{A}{A_1 & \ldots & A_n}
\]
where $A$, $A_1$, \dots, $A_n \in \Frm(\LA)$
\end{enumerate}
A formula~$F$ is a {\em theorem} of \AL{} if there is a derivation of $F$
in \AK, i.e., a finite sequence \[ F_1, F_2, \ldots, F_s = F \]
of formulas s.t.{} for each $F_i$ either
\begin{enumerate}
\item $F_i$ is a substitution instance of an axiom in $A(\AK)$, or
\item there are $F_{k_1}$, \dots,~$F_{k_n}$ with $k_j < i$
and a rule $r \in R(\AK)$, s.t.{} $F_{k_j}$ is a substitution
instance of the $j$-th premise of~$r$, and $F_i$ is a
substitution instance of the conclusion. 
\end{enumerate}
If $F$ is a theorem of~\AK{} we write $\AK \vdash F$.
The set of theorems of~\AK{} is denoted by~$\Thm(\AK)$.
\end{defn}

\begin{rem}
The above notion of a propositional rule is the one usually
used in axiomatizations of propositional logic.  It is, however,
by no means the only possible notion.  For instance, Sch\"utte's
rules
\[
\infer{A(X)}{A(\top) & A(\bot)} \qquad
\infer{A(C) \leftrightarrow A(D)}{C \leftrightarrow D}
\]
where $X$ is a propositional variable, and $A$, $C$, and $D$ are formulas,
does not fit under the above definition.
\end{rem}

\begin{defn}\label{defn:analytic}
A propositional calculus is called {\em analytic}
iff for every rule
\[
\infer[r]{A}{A_1 \ldots A_n}
\]
it holds that $\Var(A_1) \subseteq \Var(A)$, \dots,~$\Var(A_n) \subseteq \Var(A)$.
\end{defn}

\begin{rem}
Note that analytic calculi here need {\em not} have a strict
subformula property, in contrast to the notion in
sequent calculus.  Cut-free sequent calculi can easily be
be encoded in analytic Hilbert-type calculi.
Henceforth, whenever we refer to a
sequent calculus we always mean its encoding according
to the following construction.
\begin{enumerate}
\item Sequences of formulas can be coded using
a binary operator~$\cdot$.  The sequent
arrow can simply be coded as a binary operator~$\to$.
We have the following rules, to assure associativity of~$\cdot$:
\[
\infer{X \cdot \Bigl(\bigl((U\cdot V)\cdot W\bigr)\cdot Y\Bigr) \to Z}
      {X \cdot \Bigl(\bigl(U\cdot (V\cdot W)\bigr)\cdot Y\Bigr) \to Z}
\qquad
\infer{\Bigl(X \cdot \bigl((U\cdot V)\cdot W\bigr)\Bigr)\cdot Y \to Z}
      {\Bigl(X \cdot \bigl(U\cdot (V\cdot W)\bigr)\Bigr)\cdot Y \to Z}
\]
as well as the respective rules without $X$, without $Y$, without
both $X$~and~$Y$, with the rules upside-down, and also for the
right side of the sequent (20~rules total).
\item To avoid logical rules acting on sequences
instead of formulas, a formula marker~$^F$ is introduced.
Logical axioms then take the form $X^F \to X^F$.
\item The usual sequent rules can be coded using the above
constructions.
\item Some sequent rules require restrictions on the
form of the side formulas in a rule, e.g., the
R! rule in classical linear logic:
\[
\infer[$R!$]{!\Pi \to !A, ?\Gamma}{!\Pi \to A, ?\Gamma}
\]
We introduce operators $^!$ and $^?$ s.t.
\begin{enumerate}
\item $A^!$ and $B^?$ can be introduced only on $A\equiv !C$ and $B\equiv ?D$,
respectively;
\item $^!$ and $^?$ distribute over $\cdot$; and
\item $^!$ and $^?$ can always be canceled.
\end{enumerate}
R! would then take the form
\[
\infer{X^! \to {!A} \cdot Y^?}{X^! \to A \cdot Y^?}
\]
\end{enumerate}
It is easily seen that the resulting Hilbert calculus is
analytic in the sense of Definition~\ref{defn:analytic}
if the original sequent calculus was.  This also shows that
this notion of analyticity does not entail decidability,
since for instance cut-free propositional linear logic~\LL{} can be
coded in an analytic Hilbert calculus.  \LL, however, is
undecidable~\cite{LMSS:90}.
\end{rem}

\begin{ex}
Intuitionistic propositional logic is axiomatized by the following
calculus~\IPC:
\begin{enumerate}
\item Axioms:
\[
\begin{array}{ll}
a_1 & A \impl A \land A \\
a_2 & A \land B \impl B \land A \\
a_3 & A \impl B \impl (A \land C \impl B \land C) \\
a_4 & (A \impl B) \land (B \impl C) \impl (A \impl C) \\
a_5 & B \impl (A \impl B) \\
a_6 & A \land (A \impl B) \impl B \\
a_7 & A \impl A \lor B \\
a_8 & A \lor B \impl B \lor A \\
a_9 & (A \impl C) \land (B \impl C) \impl (A \lor B \impl C) \\
a_{10} & \neg A \impl A \impl B \\
a_{11} & (A \impl B) \land (A \impl \neg B) \impl \neg A \\
a_{12} & A \impl (B \impl A \land B)
\end{array}
\]
\item Rules (in usual notation):
\[
\infer[$MP$]{B}{A & A \impl B}
\]
\end{enumerate}
Gentzen's sequent calculus {\bf LJ} without cut gives an
analytical axiomatization.
\end{ex}

\begin{defn}
A {\em propositional logic}~\AL{} in the language~$\cal L$ is a
subset of~$\Frm({\cal L})$ closed under substitution.
\end{defn}

Every propositional calculus~\AK{} defines a propositional logic,
namely $\Thm(\AK) \subseteq \Frm({\cal L})$, since
$\Thm(\AK)$ is closed under substitution.
Not every propositional logic, however, is axiomatizable, let
alone finitely axiomatizable by a Hilbert calculus.
For instance, the logic
\begin{eqnarray*}
\{\Box^k(\top) & \mid & \hbox{$k$ is the G\"odel number of a} \\
& & \hbox{true sentence of arithmetic} \}
\end{eqnarray*}
is not axiomatizable, whereas the logic
\[
\{\Box^k(\top) \mid \hbox{$k$ is prime}\}
\]
is certainly axiomatizable (it is even decidable), but not by
a Hilbert calculus using only $\Box$ and $\top$. (It is easily
seen that any Hilbert calculus for $\Box$ and $\top$ has
either only a finite number of theorems or yields arithmetic
progressions of~$\Box$'s.)

\begin{defn}
A {\em propositional finite-val\-ued logic \M} is
given by a set of truth values $V(\M) = \{1$, $2$, \dots, $m\}$, the set of
{\em designated truth values} $V^+(\M) \subseteq V(\M)$, and a set of truth functions
$\tbox_j \colon V(\M)^{n_j} \to V(\M)$ for all connectives
$\Box_j \in \LA$ with arity~$n_j$.
\end{defn}

The corresponding subset of $\Frm({\cal L})$ of true formulas
is the set of tautologies of~\M, defined as follows.

\begin{defn}
A {\em valuation} \I{} is a mapping from the set of propositional
variables into~$V(\M)$.  A valuation \I{} can be extended in the standard
way to a function from formulas to truth values.  \I{} {\em
satisfies} a formula $F$, in symbols:  $\I \models_\M F$, if $\I(F)
\in V^+(\M)$.  In that case, \I{} is called a {\em model} of~$F$,
otherwise a {\em countermodel}.  A formula $F$ is a {\em tautology} of
\M{} iff it is satisfied by every valuation. Then we write $\M \models F$.
We denote the set of tautologies of \M{} by $\Taut(\M)$.
\end{defn}

\begin{ex}
The sequence
of $m$-valued G\"odel logics~\Gm{} is given by $V(\Gm) = \{0, 1, \ldots, m-1\}$,
the designated values $V^+(\Gm) = \{0\}$,
and the following truth functions:
\begin{eqnarray*}
\widetilde{\neg}_\Gm(v) & = & \cases{0 & for $v = m-1$ \cr m-1 & for $v \neq m-1$} \\
\widetilde{\lor}_\Gm(v, w) & = & \min(a, b) \\
\widetilde{\land}_\Gm(v, w) & = & \max(a, b) \\
\widetilde{\impl}_\Gm(v, w) & =& \cases{0 & for $v \ge w$ \cr w & for $v < w$}
\end{eqnarray*}
This sequence of logics was used in \cite{Godel:32} to show that
intuitionistic logic cannot be characterized by a finite matrix.
\end{ex}

In the remaining sections, we will concentrate on the
relations between calculi~\AK, logics~\AL,  and
many-valued logics~\M.  The objective is to find many-valued
logics~\M{} (or sequences thereof) that, in a sense, approximate
the calculus~\AK{} and/or the logic~\AL.

The following well-known product construction is useful for characterizing
the ``intersection'' of many-valued logics.

\begin{defn}
Let $\M$ and $\M'$ be $m$ and $m'$-valued logics, respectively.  Then
$\M \times \M'$ is the $mm'$-valued logic where $V(\M \times \M') =
V(\M) \times V(\M')$, $V^+(\M \times \M') = V^+(\M) \times V^+(\M')$,
and truth functions are defined component-wise.  I.e., if $\Box$ is an
$n$-ary connective, then \[\tbox_{\M \times \M'}(w_1,
\ldots, w_n) = \langle \tbox_\M,
\tbox_{\M'}\rangle.\]
\end{defn}

For convenience, we define the following:  Let \I{} and $\I'$ be
valuations of \M{} and $\M'$, respectively.  $\I \times \I'$ is the
valuation of $\M \times \M'$ defined by:  $(\I \times \I')(X) =
\langle \I(X), \I'(X) \rangle$.  If $\I^\times$ is a valuation of $\M
\times \M'$, then the valuations $\pi_1\I^\times$ and $\pi_2\I^\times$
of $\M$ and $\M'$, respectively, are defined by $\pi_1\I^\times(X) =
v$ and $\pi_2\I^\times(X) = v'$ iff $\I^\times(X) = \langle v,
v'\rangle$.

\begin{lem}\label{lem:tauttimes}
$\Taut(\M \times \M') = \Taut(\M) \cap \Taut(\M')$
\end{lem}

\begin{pf}
Let $A$ be a tautology of $\M \times \M'$ and $\I$ and $\I'$ be
valuations of $\M$ and $\M'$, respectively.  Since $\I \times \I'
\models_{\M \times \M'} A$, we have $\I \models_\M A$ and $\I'
\models_{\M'} A$ by the definition of~$\times$.  Conversely, let $A$
be a tautology of both \M{} and $\M'$, and let $\I^\times$ be a
valuation of $\M \times \M'$.  Since $\pi_1\I^\times \models_\M A$ and
$\pi_2\I^\times \models_{\M'} A$, it follows that $\I^\times
\models_{\M \times \M'} A$.
\end{pf}

The definition and lemma are easily generalized to the
case of finite products $\prod_i \M_i$ by induction.

When looking for a logic with as small a number of truth values
as possible which falsifies a given formula we can
use the following construction.

\begin{prop}\label{prop:red}
Let \M{} be any many-valued logic,
and $A_1$, \dots,~$A_n$ be formulas not valid in~\M.
Then there is a finite-valued logic $\M' = \Phi(\M, A_1, \ldots, A_n)$ s.t.{}
\begin{enumerate}
\item $A_1$, \dots, $A_n$ are not valid in~$\M'$,
\item $\Taut(\M) \subseteq \Taut(\M')$, and
\item $\left| V(\M') \right| \le \xi(A_1, \ldots, A_n)$, where
$\xi(A_1, \ldots, A_n) = \prod_{i=1}^n \xi(A_i)$ and $\xi(A_i)$ is the 
number of subformulas of $A_i\ + 1$. 
\end{enumerate}
This holds also if \M{} has infinitely many truth values,
provided $V(\M)$, $V^+(\M)$ and the truth functions are
recursive.
\end{prop}

\begin{pf}
We first prove the proposition for $n = 1$.
Let \I~be the interpretation in~\M{} making
$A_1$~false, and let $B_1$, \dots,~$B_r$ ($\xi(A_1) = r+1$)
be all subformulas of~$A_1$.
Every $B_i$ has a truth value~$t_i$ in~\I.
Let $\M'$ be as follows: $V(\M') = \{t_1, \ldots, t_r, \top\}$,
$V^+(\M') = V^+(\M) \cap V(\M') \cup \{\top\}$.
If $\Box \in \LA$, define $\tbox$ by
\[
\tbox(v_1, \ldots, v_n) = \cases{
t_i & if $B_i \equiv \Box(B_{j_1}, \ldots, B_{j_n})$ \cr
    & and $v_1 = t_{j_1}$, \dots, $v_n = t_{j_n}$ \cr
\top & otherwise}
\]

(1)~Since $t_r$ was undesignated in $\M$, it is also
undesignated in $\M'$. But $\I$ is also
a truth value assignment in $\M'$,
hence $\M' \not\models A_1$.

(2)~Let $C$ be a tautology of \M, and let \J~be an interpretation
in~$\M'$.  If no subformula of $C$ evaluates to $\top$
under~\J, then \J{} is also an interpretation in~\M, and
$C$ takes the same truth value in $\M'$ as in \M{} w.r.t.~\J,
which is designated also in $\M'$. Otherwise, $C$ evaluates to $\top$,
which is designated in~$\M'$. So $C$ is a tautology in~$\M'$.

(3)~Obvious.

\noindent For $n > 1$, the proposition follows by taking 
$\Phi(\M, A_1, \ldots, A_n) = \prod_{i =1}^n \Phi(\M, A_i)$
\end{pf}

Algebraic constructions can be used for simplifications of
many-valued logics.
For example, a many-valued logic~\M{} has the same tautologies
as a homomorphic image~$\M'$, if the induced congruence~$C$
on $V(\M)$ satisfies the following condition:
\[
{\rm if\quad} U \in C \quad{\rm then}\quad
V^+(\M) \cap U = \emptyset \quad{\rm or}\quad V^+(\M)
\cap (V(\M)\setminus U) = \emptyset.
\]

\section{Many-valued Covers for Calculi}

We are looking for many-valued logics~\M{}
s.t.{} $\Thm(\AK) \subseteq \Taut(\M)$.
\M{} must, however, behave ``normally''
with respect to~\AK, i.e., \AK{} must remain sound
whenever we add new operators and their truth tables to~\M{}
or add tautologies as axioms to~\AK.

\begin{defn}\label{prop:approxtest}
An $m$-valued logic~\M{} is {\em normal} for a calculus \AK{} (and \AK{}
{\em strongly sound} for~\M) if
\begin{enumerate}
\item[($*$)] All axioms $A \in A(\AK)$ are tautologies of~\M, and
for every rule~$r \in R(\AK)$: if a valuation
satisfies the premises of~$r$, it also satisfies the
conclusion.
\end{enumerate}
\M{} is then called a {\em cover} for \AK.
\end{defn}

We would like to stress the distinction between strong soundness,
a.k.a.{} normality, and soundness.  The latter is the
familiar property of a calculus to produce only valid formulas
as theorems.  This ``plain'' soundness is what we actually would
like to investigate in terms of approximations.  More precisely,
when looking for a finite-valued logic that approximates 
a given calculus, we are content if we find a logic for which
\AK{} is sound.  It is, however, not possible in general to test
if a calculus is sound for a given finite-valued logic.  It {\em is}
possible to test if it is strongly sound.  For this pragmatic reason
we consider only normal matrices for the given calculi.  The next
proposition characterizes the normal matrices in terms of strong
soundness conditions.  These are reasonable conditions which one
expects to hold of a ``normal'' matrix.

\begin{prop}
\AK{} is strongly sound for a many-valued logic~\M{}
iff $\Thm(\AK') \subseteq \Taut(\M')$ for all $\M'$ and $\AK'$, where
\begin{enumerate}
\item $\M'$ is obtained from \M{} by adding truth tables for new
operations, and
\item $\AK'$ is obtained from $\AK$ by adding tautologies of $\M'$ to
as axioms.
\end{enumerate}
\end{prop}

\begin{pf}
If:
First of all, \AK{} is sound for~\M:
Let $\AK \vdash F$.  We show that $\M \models F$ by induction on the
length~$l$ of the derivation in~\AK:

$l=1$: This means $F$ is a substi!tution instance of an axiom~$A$.

$l > 1$.  $F$ is the conclusion of a rule $r \in R(\AK)$.  If $r$ is
\[\infer[r]{A}{A_1 & \ldots & A_k}\]
and $X_1$, $X_2$, \dots, $X_n$
are all the variables in $A$, $A_1$, \dots, $A_k$, then the inference
has the form
\[
\infer{F = A[B_1/X_1, \ldots, B_n/X_n]}{A_1[B_1/X_1,\ldots, B_n/X_n]  &
\ldots & A_k[B_1/X_1,\ldots, B_n/X_n]}
\]
Let \I{} be a valuation of the variables in $F$, and let $v_i =
\I(B_i)$ ($1 \le i \le n$).  By induction hypothesis, the premises of
$r$ are valid.  This implies that, for $1 \le i \le k$, we have
$\{X_1 \mapsto v_1, \ldots, X_n \mapsto v_n\} \models A_i$.  By
hypothesis then, $\{X_1 \mapsto v_1, \ldots, X_n \mapsto v_n\} \models
A$.  But this means that $\I \models F$.  Hence, $\M \models F$.

Moreover, \AK{} satisfies conditions (1) and (2) above.

Only if:
Every axiom is derivable in \AK. By soundness, it is a tautology of~\M,
which is just what ($*$) says.  Now let $r \in R(\AK)$ be a rule, let
\I~be an interpretation which makes the premises $A_1$, \dots,~$A_k$
of $r$ true, and let $A$ be the conclusion of~$r$.
\I{} assigns truth values $v_1$, \dots,~$v_l$ to the
variables $X_1$, \dots,~$X_l$ in~$r$. Let $\M'$ be the $m$-valued logic
resulting from $\M$ by extending the language by the constants
$V_1$, \dots,~$V_l$ with values $v_1$, \dots,~$v_l$, respectively.
Let $\sigma$ be the substitution mapping $X_i$ to $V_i$.
The formulas $A_1\sigma$, \dots,~$A_l\sigma$ and (by~$r$ also) $A\sigma$
are derivable in the extension~$\AK'$ of $\AK$ by the axioms
$A_1\sigma$, \dots,~$A_l\sigma$. By (1)~and~(2), $\AK'$ is sound,
so $A\sigma$ is a tautology in $\M'$. Consequently, $\I \models A$
in~\M.
\end{pf}

\begin{cor}
If \AK~is strongly sound for~\M{} and $r$~is a directly dependent
rule of~\AK{} (i.e., $r$ can be simulated by the rules
of~\AK) then $\AK + r$ is also strongly sound for~\M.
\end{cor}

\begin{prop}\label{cor:dec}
It is decidable if a given propositional calculus is
strongly sound for a given
$m$-valued logic.
\end{prop}

Note also that for usual calculi, Property~($*$) is
relatively easy to check. For instance, modus ponens
is strongly sound iff, whenever $A$ is true,
$A \impl B$ is true iff $B$ is true;
necessitation is strongly sound if $\Box X$ is true
whenever $X$ is true.

\begin{ex}\label{ex:goedelsound}
The \IPC{} is strongly sound for the $m$-valued G\"odel logics~\Gm{}.
For instance, take axiom $a_3$:  $B \impl A \impl B$.  This is a
tautology in \Gm, for assume we assign some truth values $a$ and $b$
to $A$ and $B$, respectively.  We have two cases:  If $a \le b$, then
$(A \impl B)$ takes the value~$m-1$.  Whatever $b$ is, it certainly is
$\le m-1$, hence $B \impl A \impl B$ takes the designated value $m-1$.
Otherwise, $A \impl B$ takes the value $b$, and again (since $b \le
b$), $B \impl A \impl B$ takes the value $m-1$.

Modus ponens passes the test: Assume $A$ and $A \impl B$ both
take the value $m-1$. This means that $a \le b$. But $a = m-1$, hence
$b = m-1$.

Now consider the following extension $\G_m^\top$ of \Gm:
$V(\G_m^\top) = V(\Gm) \cup \{\top\}$, $V^+(\G_m^\top) = \{m -1, \top\}$,
and the truth functions are given by:
\[
\tbox_{\G_m^\top}(\bar{v}) = \cases{\top & if $\top \in \bar{v}$ \cr
\tbox_{\Gm}(\bar{v}) & otherwise}
\] 
for $\Box \in \{\neg, \impl, \land, \lor\}$.
Neither \IPC{} nor {\bf LJ} are strongly sound for $\G_m^\top$, but
{\bf LJ} without cut is.
\end{ex}

\begin{ex}
Consider the following calculus~{\bf K}:
\[
X \cbi \Nex X \qquad \infer[r_1]{X \cbi \Nex Y}{X \cbi Y}
\qquad \infer[r_2]{Y}{X \cbi X}
\]
It is easy to see that the corresponding logic consists
of all instances of $X \cbi \Nex^k X$ where $k \ge 1$.
This calculus is only strongly sound for the $m$-valued
logic having all formulas as its tautologies.
But if we leave out $r_2$, we can give a sequence of many-valued
logics $\M_i$, for each of which {\bf K} is strongly sound:
Take for $V(\M_n) = \{0, \ldots, n-1\}$, $V^+(\M_n) = \{0\}$,
with the following truth functions:
\begin{eqnarray*}
\widetilde{\Nex}v & = & \cases{v+1 & if $v < n-1$ \cr
n-1 & otherwise} \\
v \widetilde{\cbi} w & = & \cases{0 & if $v < w$ or $v = n-1$\cr
1 & otherwise} \\
\end{eqnarray*}
Obviously, $\M_n$ is a cover for {\bf K}.
On the other hand, $\Taut(\M_n) \neq \Frm(\LA)$, e.g., any formula
of the form $\Nex(A)$ takes a (non-designated) value~$>0$ (for $n > 1$).
In fact, every formula of the form $\Nex^k X \cbi X$
is falsified in some~$\M_n$. 
\end{ex}

\section{Optimal Covers}

By Proposition~\ref{cor:dec} it is decidable if a given $m$-valued logic~\M{}
is a cover of~\AK.  Since we can enumerate all $m$-valued logics, we can 
also find all covers of~\AK.
Moreover, comparing two many-valued logics as to their
sets of tautologies is decidable, as the next theorem will show.
Using this result, we see that we can always generate optimal
covers for \AL{}.

\begin{defn}
For two many-valued logics $\M_1$ and $\M_2$, we write
$\M_1 \bettereq \M_2$ iff $\Taut(\M_1) \subseteq \Taut(\M_2)$.

$\M_1$ is {\em better} than~$\M_2$, $\M_1 \better
\M_2$, iff $\M_1 \bettereq \M_2$ and $\Taut(\M_1) \neq \Taut(\M_2)$.
\end{defn}

\begin{thm}
Let two logics $\M_1$ and $\M_2$, $m_1$-valued and $m_2$-valued respectively,
be given. It is decidable whether $\M_1 \better \M_2$.
\end{thm}

\begin{pf}
It suffices to show the decidability of the following property:  There
is a formula~$A$, s.t.\ (*)~$\M_2 \models A$ but $\M_1 \not\models A$.
If this is the case, write $\M_1 \better^* \M_2$.
$\M_1 \better \M_2$ iff $\M_1 \better^* \M_2$ and not
$\M_2 \better^* \M_1$.

We show this by giving an upper bound on the depth of a minimal
formula~$A$ satisfying the above property.  Since the set of formulas
of~${\cal L}$ is enumerable, bounded search will produce such a
formula iff it exists.  Note that the property~(*) is decidable by
enumerating all assignments. In the following, let $m = \max(m_1,m_2)$.

Let~$A$ be a formula that satisfies~(*), i.e., there is a valuation
\I{} s.t.\ $\I {\not\models}_{\M_1} A$.  W.l.o.g.\ we can assume that
$A$~contains at most~$m$ different variables:  if it contained more,
some of them must be evaluated to the same truth value in the
counterexample~\I{} for $\M_1 \not\models A$.  Unifying these
variables leaves~(*) intact.

Let $B = \{B_1, B_2, \ldots \}$ be the set of all subformulas of~$A$.
Every formula~$B_j$ defines an $m$-valued truth function~$f(B_j)$ of
$m$~variables where the values of the variables which actually occur
in~$B_j$ determine the value of~$f(B_j)$ via the matrix of~$\M_2$.  On
the other hand, every~$B_j$ evaluates to a single truth value~$t(B_j)$
in the countermodel~\I.

Consider the formula~$A'$ constructed from $A$ as follows:  Let~$B_i$
be a subformula of~$A$ and $B_j$~be a proper subformula of~$B_i$ (and
hence, a proper subformula of~$A$).  If $f(B_i) = f(B_j)$ and $t(B_i)
= t(B_j)$, replace $B_i$ in $A$ with $B_j$.  $A'$ is shorter than~$A$,
and it still satisfies~(*).  By iterating this construction
until no two subformulas have the desired property we obtain a
formula~$A^*$.  This procedure terminates, since $A'$ is shorter than
$A$; it preserves~(*), since $A'$ remains a tautology
under~$\M_2$ (we replace subformulas behaving in exactly the same way
under all valuations) and the countermodel~\I{} is also a countermodel
for~$A'$.

The depth of $A^*$ is bounded above by $m^{m^m+1}-1$.  This is seen as
follows:  If the depth of $A^*$ is~$d$, then there is a sequence $A^*
= B_0', B_1', \ldots, B_d'$ of subformulas of $A^*$ where $B_k'$ is an
immediate subformula of $B_{k-1}'$.  Every such $B_k'$ defines a truth
function~$f(B_k')$ of $m$ variables in $\M_2$ and a truth valued
$t(B_k')$ in $\M_1$ via \I.  There are $m^{m^m}$ $m$-ary truth
functions of $m$ truth values.  The number of distinct truth
function-truth value pairs then is $m^{m^m+1}$.  If $d \ge m^{m^m+1}$,
then two of the $B_k'$, say $B_i'$ and $B_j'$ where $B_j'$ is a
subformula of $B_i'$ define the same truth function and the same truth
value.  But then $B_i'$ could be replaced by $B_j'$, contradicting the
way $A^*$ is defined.
\end{pf}

\begin{cor}
It is decidable if two many-valued logics define the same
set of tautologies. The relation $\bettereq$ is decidable.
\end{cor}

\begin{pf}
$\Taut(\M_1) = \Taut(\M_2)$ iff neither $\M_1 \better^* \M_2$ nor
$\M_2 \better^* \M_1$.
\end{pf}

Let $\simeq$ be the equivalence relation on $m$-valued logics defined
by:  $\M_1 \simeq \M_2$ iff $\Taut(M_1) = \Taut(M_2)$, and let $\MVL_m$
be the set of all $m$-valued logics over~{\cal L}.  By~${\cal M}_m$ we
denote the set of all sets~$\Taut(\M)$ of tautologies of $m$-valued
logics~$\M$.  The partial order $\langle {\cal M}_m, \subseteq\rangle$ is
isomorphic to~$\langle \MVL_m/\simeq, \bettereq/\simeq\rangle$.

\begin{prop}
$\langle {\cal M}_m, \subseteq\rangle$ is a finite complete
partial order.
\end{prop}

\begin{pf}
The set of $m$-valued logics $\MVL_m$ is obviously finite,
since there are at most $m^{n_1}m^{n_2}\cdots m^{n_c}$ different
$m$-valued matrices for~$C$. $\better$ is a partial order
on $\MVL_m/\simeq$ with the smallest element $\bot := \Frm({\cal L})$
and the largest element $\top := \emptyset$.
\end{pf}

The ``best'' logic is the one without theorems, generated by a matrix
where no connective takes a designated truth value {\em anywhere}.
The ``worst'' logic is the one where every formula of {\cal L} is a
tautology, it is generated by a matrix where every connective takes a
designated truth value {\em everywhere}.

In every complete partial order over a finite set, there exist lub and
glb for every two elements of the set.  Hence, $\langle \M, \best,
\worst, \bot, \top\rangle$ is a finite complete lattice, where
$\best$~is the lub in~$\bettereq$, and $\worst$~is the glb
in~$\bettereq$.  Since $\bettereq$~is decidable and~\M\ can be
automatically generated the functions $\best$~and~$\worst$ are
computable.

\begin{prop}
The optimal (i.e., minimal under
$\better$) $m$-valued covers of \AK{} are computable.
\end{prop}

\begin{pf}
Consider the set~$C(\AK)$
of $m$-valued covers of~\AK. Since $C(\AK)$ is finite
and partially ordered by~$\bettereq$, $C(\AK)$ contains minimal
elements. The relation $\bettereq$ is decidable, hence the
minimal covers can be computed.
\end{pf}

\begin{ex}
By Example~\ref{ex:goedelsound},
\IPC{} is strongly sound for $\G_3$.
The best 3-valued approximation of~\IPC{}
is the 3-valued G\"odel logic. In fact, it is the only 3-valued
approximation of {\em any} sound calculus~\AK{} (containing modus ponens)
for \IPL{} which has less tautologies than~\CL.
This can be seen as follows: Consider the fragment containing
$\bot$ and $\impl$ ($\neg B$ is usually defined as $B \impl \bot$).
Let \M{} be some 3-valued strongly sound approximation of~\AK.
By G\"odel's double-negation translation, 
$B$ is a classical tautology iff $\neg\neg B$ is true intuitionistically.
Hence, whenever $\M \models \neg \neg X \impl X$, then
$\Taut(\M) \supseteq \CL$. Let $0$ denote the value of $\bot$ in \M,
and let $1 \in V^+(\M)$. We distinguish cases:
\begin{enumerate}
\item $0 \in V^+(\M)$: Then $\Taut(\M) = \Frm(\LA)$, since $\bot \impl X$
is true intuitionistically, and by modus ponens: $\bot, \bot \impl X / X$.
\item $0 \notin v^+(\M)$: Let $u$ be the third truth value.
\begin{enumerate}
\item $u \in V^+(\M)$: Consider $A \equiv ((X \impl \bot) \impl \bot) \impl X$.
If $\I(X)$ is $u$~or~$1$, then, since everything implies something true,
$A$ is true (Note that we have $Y, Y \impl (X \impl Y) / X \impl Y$).
If $\I(X) = 0$, then (since $0 \impl 0$ is true, but $u \impl 0$ and $1 \impl 0$
are both false), $A$ is true as well. So $\Taut(\M) \supseteq \CL$.
\item $u \notin V^+(\M)$, i.e., $V^+(\M) = \{1\}$: Consider the
truth table for implication.  Since $B \impl B$, $\bot \impl B$,
and something true is implied by everything, the upper right triangle
is~$1$. We have the following table:
\[
\begin{array}{c|ccc}
\impl & 0 & u & 1 \\ \hline
0 &     1 & 1 & 1 \\
u &     v_1 & 1 & 1 \\
1 &     v_0 & v_2 & 1
\end{array} 
\]
Clearly, $v_0$ cannot be~$1$. If $v_0 = u$, we have, by $((X \impl X) \impl \bot) \impl Y$,
that $v_1 = 1$. In this case, $\M \models A$ and hence $\Taut(\M) \supseteq \CL$.
So assume $v_0 = 0$.
\begin{enumerate}
\item $v_1 = 1$: $\M \models A$ (Note that only the case of $((u \impl 0) \impl 0) \impl u$
has to be checked).
\item $v_1 = u$: $\M \models A$.
\item $v_1 = 0$: With $v_2 = 0$, \M{} would be incorrect ($u \impl (1 \impl u)$ is false).
If $v_2 = 1$, again $\M \models A$. The case of $v_2 = u$ is the G\"odel logic,
where $A$ is not a tautology.  
\end{enumerate}
\end{enumerate}
\end{enumerate}
\end{ex}

Note that it is in general impossible to algorithmically construct
a $\bettereq$-minimal $m$-valued {\em logic}~\M{} (i.e., given
independently of a calculus)
with $\AL \subseteq \Taut(\M)$, because, e.g.,
it is undecidable whether \M{} is empty or not: e.g.,
take \[\AL = \cases{\{\Box^k(\top)\} & if $k$ is the least solution of $D(x) = 0
$\cr
\emptyset & otherwise}\]
where $D(x) = 0$ is the diophantine representation of some undecidable
set.

\section{Sequential Approximations of Calculi}

In the previous section we have shown that it is always possible to
obtain the best $m$-valued covers of a given calculus,
but there is no way to tell {\em how
good} these covers are.  In this section, we investigate
the relation between sequences of many-valued logics
and the set of theorems  of
a calculus~\AK.  Such sequences are called
{\em sequential approximations} of \AK{} if they
verify all theorems and refute all non-theorems of~\AK. 
Put another way, this is a question about the
limitations of Bernays' method.
On the negative side an immediate result says that calculi
for undecidable logics do not have sequential approximations.
If, however, a propositional logic is decidable, it also
has a sequential approximation (independent of a calculus).
However, they all have a uniquely defined {\em many-valued closure},
whether they are decidable or not.  This is the set of all sentences
which cannot be proved underivable using a Bernays-style many-valued argument.
If a calculus has a sequential
approximation, then the set of its theorems equals its many-valued closure.
If it does not, then its closure is a proper superset.
Different calculi for one and the same logic may have
different many-valued closures according to their degree of analyticity.

\begin{defn}\label{defn:sapprox}
Let \AK{} be a calculus and
let $\A = \langle \M_1, \M_2, \M_3, \ldots, \M_j, \ldots\rangle$ ($j \in
\omega$) be a sequence of many-valued logics s.t.
\begin{enumerate}
\item \A~is given by a recursive procedure,
\item $\M_i \bettereq \M_j$ iff $i \ge j$, and
\item $\M_i$ is a cover for \AK.
\end{enumerate}
\A{} is called a {\em sequential approximation} of \AK{} iff
$\Thm(\AK) = \bigcap_{j \in \omega} \Taut(\M_j)$.
We say \AK{} is {\em approximable}, if there is such a sequential
approximation for~\AK{}.
\end{defn}

Condition (2) above is technically not necessary.  Approximating sequences
of logics in the literature (see next example), however, satisfy this
condition.  Furthermore, with the emphasis on ``approximation,'' it seems more
natural that the sequence gets successively ``better.''

\begin{ex}
Consider the sequence $\G = \langle \G_i \rangle_{i\ge 2}$ of G\"odel
logics and intuitionistic propositional logic \IPC.
$\Taut(\G_i) \supset \Thm(\IPC)$, since $\G_i$ is a cover
for~\IPC. Furthermore, $\G_{i+1} \better \G_i$. This
has been pointed out by \cite{Godel:32}, for a detailed
proof see \cite[Satz~3.4.1]{Gottwald:89}.  It
is, however, not a sequential approximation of \IPC: The formula $(A \impl B)
\lor (B \impl A)$, while not a theorem of \IPL, is a tautology
of all $\G_i$. In fact, $\bigcap_{j \ge 2} \Taut(\G_i)$ is the
set of tautologies of the infinite-valued G\"odel logic~$\G_\aleph$,
which is axiomatized by the rules of \IPC{} plus the above formula.
This has been shown in \cite{Dummett:59} (see also \cite[\S~3.4]{Gottwald:89}).
Hence, \G{} is a sequential approximation of $\G_\aleph = \IPC +
(A \impl B) \lor (B \impl A)$.

Ja{\'s}kowski \cite{Jaskowski:36} gave a sequential approximation of \IPC.
That \IPC~is approximable is also a consequence of Theorem~\ref{thm:fmpapprox},
with the proof adapted to Kripke semantics for intuitionistic propositional
logic, since \IPL~has the
finite model property \cite[Ch.~4, Theorem~4(a)]{Gabbay:81}.
\end{ex}

The natural question to ask is: Which calculi are approximable?
First we give the unsurprising negative answer for undecidable calculi.

\begin{prop}\label{prop:undecapprox}
If \AK{} is undecidable, then it is not approximable.
\end{prop}

\begin{pf}
If \AK{} were approximable, there were a sequence $\A = \langle \M_1, \M_2, \M_3,
\ldots\rangle$ s.t. $\bigcap_{j\ge 2} \Taut(\M_j) = \Thm(\AK)$.
If $N$ is a non-theorem of~\AK{}, then there would be an
index~$i$ s.t.\ $N$ is false in $\M_i$.  But this would yield a
semi-decision procedure for non-theorems of~\AK: Try for each~$j$
whether $N$~is false in $\M_j$. If $N$ is a non-theorem, this will be established
at $j = i$, if not, we may go on forever.  This contradicts the
assumption that the non-theorems of \AK{} are not r.e.\
(\AK{} is undecidable and the theorems are r.e.).
\end{pf}

\begin{ex}
This shows that a result similar to that for~\IPC{}
cannot be obtained for full propositional linear logic.
\end{ex}

If \AK{} is not approximable (e.g., if it is undecidable),
then the intersection of all covers for~\AK{}
is a proper superset of~$\Thm(\AK)$.  This intersection
has interesting properties.

\begin{defn}
The {\em many-valued closure}~$\MC(\AK)$ of a calculus~\AK{}
is the set of formulas which are true in every many-valued
cover for~\AK{}.
\end{defn}

$\MC(\AK)$ is unique, since it obviously equals $\bigcap_{\M \in S} \Taut(\M)$
where $S$ is the set of all covers for~\AK.
It is also approximable, an approximating sequence is given by
\begin{eqnarray*}
\M_1 & = & \M'_1 \\
\M_i & = & \M_{i-1} \times \M_i'
\end{eqnarray*}
where $\M_i'$ is an enumeration of~$S$.

The many-valued closure,
however, need not be trivial (i.e., equal to $\Frm(\LA)$)---even for
undecidable~\AK{}.

\begin{prop}
If \AK{} is analytical then $\MC(\AK)$ is decidable.
\end{prop}

\begin{pf}
Assume $\AK{}$ is analytical.
A decision procedure for $A \in \MC(\AK)$ is given by the following:
Enumerate all many-valued logics~$\M_i$ in order of increasing
number of truth values.  Check if $\AK$ is strongly sound for $\M_i$
(decidable by Proposition~\ref{cor:dec}).  If it is strongly sound,
then check whether $\M_i \models A$.  If not, terminate with $A \notin \MC(\AK)$.
By Proposition~\ref{prop:red}, we only have to search until
all many-valued logics with number of truth values $\le \xi(A)$ have
been checked, provided \AK{} is strongly sound for $\M' = \Phi(\M, A)$.
Since $A$ must be a non-tautology
of some cover~$\M$ of~$\AK$ for $A \notin \MC(\AK)$ to  hold,
we can assume that $\M$ is a cover of $\AK$.
Since $\Taut(\M) \subseteq \Taut(\M')$, all axioms of~\AK{}
are tautologies in~$\M'$.  Let
\[
\infer[r]{A}{A_1 \ldots A_n}
\]
be a rule in \AK, and let \J{} be an interpretation
in $\M'$ making each $A_j$ true. If \J{} maps no variable
to $\top$, \J~is also an interpretation in $\M$.  Then, since
\AK{} is sound for \M, $A$~is true under~\J{} (in both $\M$~and~$\M'$).
Otherwise, if \J{} assigns $\top$ to some variable~$X$,
$A$~is true under~\J{} since $X$~occurs in~$A$ (recall
that \AK{} is analytical).  So \AK{} is strongly sound for~$\M'$.
\end{pf}

\begin{cor}
The many-valued closure of cut-free propositional linear logic \LL{} is decidable.
\end{cor}

\begin{cor}
If \AK{} is analytic and decidable, then $\MC(\AK) = \Thm(\AK)$.
\end{cor}

\begin{pf}
Certainly $\Thm(\AK) \subseteq \MC(\AK)$. Let $A \notin \Thm(\AK)$.  Then the
(infinite-valued) Lindenbaum logic~$\AL(\AK)$ \cite[Satz~3]{LukasiewiczTarski:30}
for \AK{} falsifies~$A$.  Since \AK{} is decidable, $\AL(\AK)$ is
effectively given.
$\AL(\AK)$ satisfies ($*$).
It is easy to see that $\Phi(\AL(\AK), A)$ also satisfies~($*$).
By Proposition~\ref{prop:red} and the argument of the above
proof, there is a finite-valued cover
for \AK{} falsifying~$A$. Hence, $A \notin \MC(\AK)$.
\end{pf}

The last corollary can be used to uniformly obtain semantics
for decidable analytic Hilbert calculi.

\section{Sequential Approximations of Other Representations}

Propositional logic can also be given by effective representations
other than calculi.  A decidable logic, for instance, may be represented
by a decision procedure.  Logics with Kripke semantics which have
the finite model property can be given by the r.e.{} sequence
of their finite models.  In this section, we investigate
the question of sequential approximation for these representations.

\begin{prop}
For every decidable propositional logic~\AL{} there is a sequence~\A{}
of many-valued logics~$\M_i$ satisfying 
\begin{enumerate}
\item \A~is given by a recursive procedure,
\item $\M_i \bettereq \M_j$ iff $i \ge j$, and
\item $\AL \subseteq \Taut(\M_i)$,
\end{enumerate}
s.t.{} $\AL = \bigcap_{i \ge 2} \Taut(\M_i)$.
\end{prop}

\begin{pf}
The proof uses an argument similar to that of
Lindenbaum \cite[Satz~3]{LukasiewiczTarski:30}.
Let $\Frm_i(\LA) \subset \Frm(\LA)$ be the set of formulas
of depth $\le i$ (which is finite up to renaming of variables).
To every formula $F \in \Frm(\LA)$ we assign a code
$\cd{F}$, yielding the sets $\cd{\Frm_i(\LA)}$ for all $i \in \omega$.
We construct a sequential approximation of~\AL{}
as follows: $V(\M_i) = \cd{\Frm_i(\LA)} \cup \{\top\}$, with the
designated values $V^+(\M_i) = \cd{\Frm_i(\LA)} \cap \cd{\AL} \cup \{\top\}$.
The truth tables for $\M_i$ are given by:
\[\begin{array}{@{}l@{}}
\tbox_{\M_i}(v_1, \ldots, v_n) = \\
\quad = \cases{\cd{\Box(F_1, \ldots, F_n)} & if $v_j = \cd{F_j}$ for $1 \le j \le n$ \cr
& and $\Box(F_1, \ldots, F_n) \in \Frm_i(\LA)$ \cr
\top & otherwise}
\end{array}\]
$\M_i$ is constructed in such a way as to agree with \AL{} on all
formulas of depth~$\le i$, and to make all formulas
of depth~$> i$ true.  Hence, $\Taut(\M_i) \supseteq \AL$, and
$\M_i \bettereq \M_{i+1}$. Every formula~$F$ false in \AL{}
is also false in some $\M_i$ (namely in all $\M_i$ with $i \ge$ the depth
of~$F$).   
\end{pf}

Note that it is in general impossible to algorithmically construct
a $\bettereq$-minimal $m$-valued logic~\M{}
with $\AL \subseteq \Taut(\M)$, because, e.g.,
it is undecidable whether \M{} is empty or not: e.g.,
take \[\AL = \cases{\{\Box^k(\top)\} & if $k$ is the least solution of $D(x) = 0$\cr
\emptyset & otherwise}\]
where $D(x) = 0$ is the diophantine representation of some undecidable
set.

The following definitions are taken from \cite{Chellas:80}.

\begin{defn}
A {\em modal logic}~\AL{} has as its language~$\cal L$
the usual propositional connectives plus two unary
{\em modal operators:} $\Box$ (necessary) and $\Diamond$ (possible).
A {\em Kripke model} for~$\cal L$ is a triple $\langle W, R, P\rangle$,
where
\begin{enumerate}
\item $W$ is any set: the set of {\em worlds},
\item $R \subseteq W^2$ is a binary relation on $W$: the {\em accessibility relation},
\item $P$ is a mapping from the propositional variables to subsets of~$W$.
\end{enumerate}
A modal logic \AL{} is characterized by a class of Kripke models for \AL.
\end{defn}

This is called the {\em standard semantics} for modal logics
(see \cite[Ch.~3]{Chellas:80}).  The semantics of formulas in standard
models is defined as follows:

\begin{defn}
Let \AL{} be a modal logic, $\K_\AL$ be its characterizing class of Kripke
models. Let $K = \langle W, R, P\rangle \in \K_\AL$ be a Kripke model
and $A$ be a modal formula.

If $\alpha \in W$ is a possible world, then we say $A$ is {\em true in $\alpha$},
$\alpha \models_\AL A$, iff the following holds:
\begin{enumerate}
\item $A$ is a variable: $\alpha \in P(X)$
\item $A \equiv \neg B$: not $\alpha \models_\AL B$
\item $A \equiv B \land C$: $\alpha \models_\AL B$ and $\alpha \models_\AL C$
\item $A \equiv B \lor C$: $\alpha \models_\AL B$ or $\alpha \models_\AL C$
\item $A \equiv \Box B$: for all $\beta \in W$ s.t. $\alpha \mathrel{R} \beta$ it holds that $\beta \models_\AL B$
\item $A \equiv \Diamond B$: there is a $\beta \in W$ s.t. $\alpha \mathrel{R} \beta$ and $\beta \models_\AL B$
\end{enumerate}
We say $A$ is {\em true} in $K$, $K \models_\AL A$, iff for all $\alpha \in W$
we have $\alpha \models_\AL A$. $A$ is {\em valid in \AL}, $\AL \models A$,
iff $A$ is true in every Kripke model $K \in \K_\AL$. By $\Taut(\AL)$ we
denote the set of all formulas valid in \AL.
\end{defn}

Many of the modal logics in the literature have
the {\em finite model property (fmp)}: for every $A$ s.t.\
$\AL \not\models A$, there is a
finite Kripke model $K = \langle W, R, P\rangle \in \K$ (i.e., $W$ is finite),
s.t.\ $K \not\models_\AL A$ (where \AL{} is characterized by~$\K$).
We would like to exploit the fmp  to construct 
sequential approximations. This can be done as follows:

\begin{defn}\label{defn:MA}
Let $K = \langle W, R, P\rangle$ be an effectively given
finite Kripke model.  We define
the many-valued logic $\M_K$ as follows:

\begin{enumerate}

\item $V(\M_K) = \{0, 1\}^W$, the set of 0-1-sequences with indices
from $W$.

\item $V^+(\M_K) = \{1\}^W$, the singleton of the sequence constantly
equal to~1.

\item $\widetilde\neg_{\M_K}$, $\widetilde\lor_{\M_K}$,
$\widetilde\land_{\M_K}$, $\widetilde\impl_{\M_K}$ are defined
componentwise from the classical truth functions

\item $\tbox_{\M_K}$ is defined as follows:
\[
\tbox_{\M_K}(\langle w_\alpha \rangle_{\alpha \in W})_\beta =
\cases{1 & if for all $\gamma$ s.t.\cr & $\beta \mathrel{R} \gamma$, $w_\gamma = 1$ \cr
0 & otherwise}
\]
\item $\tdiamond_{\M_K}$ is defined as follows:
\[
\tdiamond_{\M_K}(\langle w_\alpha \rangle_{\alpha \in W})_\beta =
\cases{1 & if there is a $\gamma$ s.t.\cr
& $\beta \mathrel{R} \gamma$ and $w_\gamma = 1$ \cr
0 & otherwise}
\]
\end{enumerate}
Furthermore, $\I_K$ is the valuation defined by
$\I_K(X)_\alpha = 1$ iff $\alpha \in P(X)$ and $= 0$ otherwise.
\end{defn}

\begin{lem}\label{lem:MKmodels}
Let $\AL$ and $K$ be as in Definition~\ref{defn:MA}.
Then the following hold:
\begin{enumerate}
\item Every valid formula of \AL{} is a tautology of $\M_K$.
\item If $K \not\models_\AL A$ then $\I_K \not\models_{\M_K} A$.
\end{enumerate}
\end{lem}

\begin{pf}
Let $B$ be a modal formula, and $K' = \langle W, R, P'\rangle$.
We prove by induction that $\val_{\I_{K'}}(B)_\alpha = 1$ iff
$\K' \models_\AL B$:

$B$ is a variable: $P'(B) = W$ iff $\I_K(B)_\alpha = 1$ for all
$\alpha \in W$ by definition of~$\I_K$.

$B \equiv \neg C$: By the definition of $\widetilde\neg_{\M_K}$,
$\val_{\I_K}(B)_\alpha = 1$ iff $\val_{\I_K}(C)_\alpha = 0$.
By induction hypothesis, this is the case iff $\alpha \not\models_\AL C$.
This in turn is equivalent to $\alpha \models_\K B$. Similarly
if $B$ is of the form $C \land D$, $C \lor D$, and $C \impl D$.

$B \equiv \Box C$: $\val_{\I_K}(B)_\alpha = 1$ iff for all $\beta$
with $\alpha \mathrel{R} \beta$ we have $\val_{\I_K}(C)_\beta = 1$.
By induction hypothesis this is equivalent to
$\beta \models_\AL C$. But by the definition of $\Box$ this
obtains iff $\alpha \models_\AL B$. Similarly for $\Diamond$.

(1)~Every valuation $\I$ of $\M_K$ defines a function $P_\I$ via
$P_\I(X) = \{\alpha \mid \I(X)_\alpha = 1\}$. Obviously,
$\I = \I_{P_\I}$.
If $\AL \models B$, then $\langle W, R, P_\I\rangle \models_\AL B$.
By the preceding argument then $\val_\I(B)_\alpha = 1$ for
all $\alpha \in W$. Hence, $B$ takes the
designated value under every valuation.

(2)~$A$ is not true in $K$. This is the case only if
there is a world $\alpha$ at which
it is not true. Consequently, $\val_{\I_K}(A)_\alpha = 0$ and $A$~takes
a non-designated truth value under~$\I_K$.
\end{pf}

The above method can be used quite in general to construct
many-valued logics from Kripke structures for not only modal
logics, but also for intuitionistic logic.
Kripke semantics for \IPL{} are defined quite similar,
with the exception that $\alpha \models A \impl B$
iff $\beta \models A \impl B$ for all $\beta \in W$ s.t.{}
$\alpha \mathrel{R} \beta$.  \IPL{} is then characterized by the
class of all finite trees \cite[Ch.~4, Thm.~4(a)]{Gabbay:81}.
Note, however, that for intuitionistic Kripke semantics
the form of the {\em assignments}~$P$ is restricted:
If $w_1 \in P(X)$ and $w_1 \mathrel{R} w_2$ then
also $w_2 \in P(X)$ \cite[Ch.~4, Def.~8]{Gabbay:81}.
Hence, the set of truth values has to be restricted in a similar way.
Usually, satisfaction for intuitionistic Kripke semantics
is defined by satisfaction in the {\em initial} world.
This means that every sequence where the first entry equals~1
should be designated.  By the above restriction,
the only such sequence is the constant 1-sequence.

\begin{ex}
The Kripke tree with three worlds
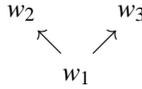
\[
\begin{array}{r@{}c@{}c@{}c@{}l}
w_2 & & & & w_3 \\
 & \nwarrow & & \nearrow \\
& & w_1
\end{array}
\]  
yields a five-valued logic~${\bf T}_3$, with
$V({\bf T_3}) = \{000, 001, 010, 011, 111\}$,
$V^+({\bf T_3}) = \{111\}$,
the truth table
for implication 
\[
\begin{array}{c|cccccccc}
\impl & 000 & 001 & 010 & 011 &  111 \\
\hline
000 & 111 & 111 & 111 & 111 &  111 \\
001 & 010 & 111 & 010 & 111 &  111 \\
010 & 001 & 001 & 111 & 111 &  111 \\
011 & 000 & 001 & 010 & 111 &  111 \\
111 & 000 & 001 & 010 & 011 &  111
\end{array}\]
$\bot$ is the constant $000$, $\neg A$ is defined
by $A \impl \bot$, and $\lor$ and $\land$ are
given by the componentwise classical operations.

The Kripke chain with four worlds corresponds
directly to the five-valued G\"odel logic~${\bf G}_5$.
It is well know that $(X \impl Y) \lor (Y \impl X)$
is a tautology in all \Gm.  Since ${\bf T}_3$
falsifies this formula (take $001$ for $X$ and
$010$ for $Y$), we know that ${\bf G}_5$
is not the best five-valued approximation of~\IPL.

Furthermore, let 
\begin{eqnarray*}
O_5 & = & \bigwedge_{1 \le i < j \le 5} (X_i \impl X_j) \lor (X_j \impl X_i) {\rm \ and}\\
F_5 & = & \bigvee_{1 \le i < j \le 5} (X_i \impl X_j).
\end{eqnarray*}
$O_5$ assures that the truth values assumed by $X_1$, \ldots, $X_5$ are
linearly ordered by implication.  Since neither $010 \impl 001$ nor $001 \impl 010$
is true, we see that there are only four truth values which can
be assigned to $X_1$, \dots,~$X_5$ making $O_5$ true.
Consequently, $O_5 \impl F_5$ is valid in ${\bf T}_3$.
On the other hand, $F_5$ is false in ${\bf G}_5$.
\end{ex}
 
\begin{thm}\label{thm:fmpapprox}
Let \AL{} be a modal logic characterized by a r.e.{} set of finite Kripke models,
and $\langle A_1, A_2, \ldots\rangle$ an enumeration
of its non-theorems. A sequential approximation of~\AL{} is given by
$\langle \M_1, \M_2, \ldots \rangle$ where
$\M_1 = \M_{K_1}$, and $\M_{i+1} = \M_i \times \M_{K_{i+1}}$
where $K_i$ is the smallest finite model s.t.\ $K_i \not\models_\AL A_i$
\end{thm}

\begin{pf}
(1)~$\Taut(\M_i) \supseteq \Taut(\AL)$:  By induction on $i$:  For $i =
1$ this is Lemma~\ref{lem:MKmodels}~(1).  For $i > 1$ the statement
follows from Lemma~\ref{lem:tauttimes}, since $\Taut(\M_{i-1})
\supseteq \Taut(\AL)$ by induction hypothesis, and $\Taut(\M_{K_i})
\supseteq \Taut(\AL)$ again by Lemma~\ref{lem:MKmodels}~(1).

(2)~$\M_i \bettereq \M_{i+1}$ from $A \cap B \subseteq A$ and
Lemma~\ref{lem:tauttimes}.

(3)~$\Taut(\AL) = \bigcap_{i\ge 1} \Taut(\M_i)$.  The
$\subseteq$-direction follows immediately from~(1).  Furthermore, by
Lemma~\ref{lem:MKmodels}~(2), no non-tautology of~\AL{} can be a member
of all $\Taut(\M_i)$, whence $\supseteq$ holds.
\end{pf}

\begin{rem}
Note that Theorem~\ref{thm:fmpapprox} does not hold in general if
\AL{} is not finitely axiomatizable.  This follows from
Proposition~\ref{prop:undecapprox} and the existence of an undecidable
recursively axiomatizable modal logic which has the fmp (see
\cite{Urquhart:81}).  Note also the condition in Theorem~\ref{thm:fmpapprox}
that there is an enumeration of the non-theorems of~\AL.  Since
finitely axiomatizable logics with the fmp are decidable (\cite{Harrop:58}),
there always is such an enumeration for the logics we consider.
\end{rem}

This theorem can also be used to show that the many-valued
closure of a calculus for a modal logic with the fmp equals
the logic itself, provided that the calculus
contains modus ponens and necessitation as
the only rules. (All standard axiomatizations are of this form.)

\section{Conclusion}

The main open problem, especially in view of possible applications
in computer science,
is the complexity of the computation of optimal covers.
One would expect that it is tractable at least for
some reasonable classes of calculi which are syntactically
characterizable, e.g., analytic calculi.

A second problem is in how far approximations can be
found for first-order logics and calculi.  One obstacle, for
instance, is that it is difficult to check whether a matrix
is normal for a given calculus, in particular if the rules
of the calculus are not ``monadic'' in the sense that they
manipulate more than one variable at a time.  In any case,
a systematic treatment only seems feasible for many-valued
logics with, at most, distribution quantifiers~\cite{Carnielli:87}.

\bibliography{logic}

\end{document}